\RequirePackage{ifpdf}
\ifpdf 
\documentclass[pdftex]{sigma}
\else
\documentclass{sigma}
\fi

\begin{document}

\renewcommand{\PaperNumber}{021}

\FirstPageHeading

\ShortArticleName{Pachner Move $3\to 3$ and Affine
Volume-Preserving Geometry in $\mathbb R^3$}

\ArticleName{Pachner Move $\boldsymbol{3\to 3}$ and Affine
Volume-Preserving\\ Geometry in $\boldsymbol{\mathbb R^3}$}

\Author{Igor G. KOREPANOV} \AuthorNameForHeading{I.G. Korepanov}

\Address{South Ural State University, 76 Lenin Ave., 454080
Chelyabinsk, Russia}
\Email{\href{mailto:kig@susu.ac.ru}{kig@susu.ac.ru}}

\ArticleDates{Received October 06, 2005, in final form November
21, 2005; Published online November 24, 2005}

\Abstract{Pachner move $3\to 3$ deals with triangulations of
four-dimensional manifolds.
 We present an algebraic relation corresponding in a natural way to this move and based,
 a~bit paradoxically, on {\em three\/}-dimensional geometry.}


\Classification{57Q99; 57M27; 57N13}

\section{Introduction}
\label{sec intro}

Pachner moves are elementary rebuildings of a manifold
triangulation. One can transform any triangulation of a fixed
manifold into any other triangulation using a sequence of enough
Pachner moves\footnote{See the excellent paper~\cite{Lickorish}
for exact theorems.}.
 In any given dimension~$n$ of a manifold, there exist a finite number, namely $n+1$,
 types of Pachner moves. In dimension~4, with which we will
deal here, these are moves $3\to 3$, $2\leftrightarrow 4$
 and $1\leftrightarrow 5$, where the numbers show how many 4-simplices
 are taken away from the triangulation and with how many simplices they are replaced.

In general, if ``elementary rebuildings'' of some object are
defined, then there is an interesting problem of finding algebraic
formulas corresponding, in a sense, to these rebuildings. For
Pachner moves, the immediate aim of such formulas may be the
construction of manifold invariants based on them. It appears,
however, that the real significance of such formulas is both wider
and deeper. A good example of this is provided by the Onsager
star-triangle relation and its generalization --- Yang--Baxter
equation. The latter proves to be relevant, in its various
versions, to as different scientific areas as statistical physics,
classical soliton equations, knot theory and the theory of
left-symmetric algebras\footnote{One can consult the
paper~\cite{leftsymmetric} about Yang--Baxter equation in the
theory of left-symmetric algebras.}.

Algebraic relations between quantities of geometric origin,
corresponding in a natural way to Pachner moves in dimensions 3
and~4, have been constructed in
papers~\cite{arx_f-la_1,arx_f-la_2,Korepanov_JNMP,SL2}. These
relations deserve the name of ``geometric-semiclassical'' because,
in the case of three-dimensional Euclidean
geometry~\cite{arx_f-la_1,Korepanov_JNMP}, they can be obtained by
a semiclassical limit from the pentagon equation for $6j$-symbols,
according to the Ponzano--Regge--Roberts
formula~\cite{Ponzano-Regge,Roberts}\footnote{See also similar
formulas for spherical geometry in~\cite{TW}.}.
 And when we are using other kinds of geometry, such as area-preserving plane geometry in~\cite{SL2,TMP2004},
 or in the four-dimensional Euclidean case~\cite{arx_f-la_2,TMP4dim1,TMP4dim2,TMP4dim3},
 we think that our formulas still deserve to be called geometric-semiclassical:
 the constructions are always very clearly akin to those in the three-dimensional Euclidean case.
 The problem is, however, that the corresponding quantum formulas are not known!
 This problem seems to be hard (its solution will perhaps lead to new topological
 field theories in the spirit of Atiyah's axioms~\cite{Atiyah}), but it is worth some effort.
 To begin, one can search for as many as possible (semi)classical formulas,
 having in mind the following way of quantization known from the experience
 in studying integrable models in mathematical physics:
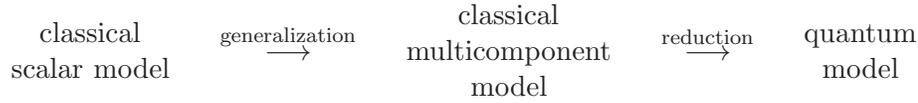
\[
\begin{array}{c}\rm classical\\ \rm scalar\ model\end{array}
\quad \stackrel{\rm generalization}{\longrightarrow} \quad
\begin{array}{c}\rm classical\cr \rm multicomponent\cr \rm model\end{array}
\quad \stackrel{\rm reduction}{\longrightarrow} \quad
\begin{array}{c}\rm quantum\cr \rm model\end{array}\;.
\]

Presenting here a ``four-dimensional'' relation, we mean not just
its applications to topology of four-manifolds but also possible
reduction(s) from four to three dimensions. In general, we have in
mind the following way of constructing large families of manifold
invariants: starting from a manifold~$M$, we construct a
simplicial complex~$K$ of a greater dimension --- the reader can
imagine a double cone over~$M$ just for a simple example --- and
consider a suitable invariant of~$K$. We hope to demonstrate in
one of future papers that this leads to interesting algebraic
structures already for a two-dimensional~$M$ (even though we know
everything about two-dimensional topology).

The direct aim of the present note is generalization of the
``SL(2)-solution of pentagon equation'' from~\cite{SL2,TMP2004}.
Here ``solution of pentagon equation'' means an algebraic relation
corresponding in a natural sense to a move $2\to 3$ of a
three-dimensional manifold triangulation; instead of group SL(2),
we now prefer to relate the constructions of~\cite{SL2,TMP2004} to
the group of affine area-preserving transformations of the
plane~$\mathbb R^2$. Below, we only generalize the pentagon
relation of paper~\cite{SL2}; we plan to present the analogue of
the acyclic complex from~\cite{TMP2004}, which yields the actual
manifold invariant, in a later publication.

\section[The 3 -> 3 relation: recalling the Euclidean case]{The $\boldsymbol{3\to 3}$
relation: recalling the Euclidean case} \label{sec_recall}

The move $3\to 3$ occupies, one can say, the central place among
the Pachner moves in dimension four. Under this move, a cluster of
three adjacent 4-simplices is replaced with a cluster of other
three adjacent 4-simplices having the same boundary. We will
construct an algebraic relation corresponding to this move,
satisfied by some quantities~$\lambda$ related to tetrahedron
volumes in a~{\em three\/}-dimensional space~$\mathbb R^3$.
Besides these 3-volumes, we will be using linear dependencies
among vectors joining points in~$\mathbb R^3$ (e.g., in
formula~(\ref{eq_c}) below). Both volumes and linear dependencies
are preserved under the action of the group of {\em
volume-preserving affine transformations}, and in this sense our
constructions are based on three-dimensional volume-preserving
affine geometry.

An analogous algebraic relation, but using {\em four-dimensional
Euclidean\/} geometry, was proposed in paper~\cite{TMP4dim1}, and
then in~\cite{TMP4dim2,TMP4dim3} an invariant of 4-manifolds has
been constructed based on this. As the analogy with Euclidean case
was actually guiding us when dealing with other geometries, it
makes sense to recall here, after some necessary explanations, the
formula~(8) from paper~\cite{TMP4dim1}.

In considering the Pachner move $3\to 3$, we will be dealing with
six vertices in the triangulation, we call them $A$, $B$, $C$,
$D$, $E$, and~$F$. The sequence of letters $ABCDEF$ with, say,
letter~$D$ taken away, will be denoted as~$\hat D$. So, the
Pachner move $3\to 3$ starts with 4-simplices $ABCDE=\hat F$,
$ABCFD=-\hat E$ and $ABCEF=\hat D$ (the minus sign is responsible
for the orientation, and we choose a consistent orientation for
4-simplices) having the common two-dimensional face~$ABC$, and
replaces them with 4-simplices $\hat A$, $-\hat B$ and~$\hat C$
having the common two-dimensional face~$DEF$ and the same common
boundary.

Imagine now that points $A,\ldots,F$ are placed in the
space~$\mathbb R^4$ {\em equipped with Euclidean metric}. This
means, in particular, that the edges of 4-simplices acquire
Euclidean lengths, the two-dimensional faces --- areas (we denote,
e.g., the area of $ABC$ as $S_{ABC}$), and 4-simplices ---
4-volumes (we denote the 4-volume of $\hat A$ as $V_{\hat
A}$)\footnote{It is interesting to note that no {\em
three\/}-volumes enter in the formula \cite[(8)]{TMP4dim1},
presented below as (\ref{eq_z}). In contrast with this, our
``affine volume-preserving'' formula~(\ref{eq_j}) is based exactly
on three-volumes.}.
 The lengths of 15 edges joining the points $A,\ldots,F$ are not independent ---
 they satisfy the well-known relation: the {\em Cayley--Menger determinant\/} made of these lengths must vanish.

We will, however, represent this relation in a different form,
namely in terms of {\em dihedral angles at a two-dimensional
face}. It is clear that Euclidean simplices $\hat F$, $-\hat E$
and~$\hat D$, if the lengths of their 15 edges are given, can be
placed together in~$\mathbb R^4$ if and only if the sum
$\vartheta_{ABC}^{\hat F} + \vartheta_{ABC}^{-\hat E} +
\vartheta_{ABC}^{\hat D}$ of dihedral angles at face~$ABC$ in
these three 4-simplices is~$2\pi$ (to be exact, their {\em
algebraic\/} sum must be $0$ modulo~$2\pi$: a dihedral angle must
be taken sometimes with a minus sign --- see
\cite[Section~3]{TMP4dim1} for details). What will happen if we
slightly, but arbitrarily, deform these 15 lengths?  Each
4-simplex individually still can be put in~$\mathbb R^4$, of
course, but if one tries to put them all together, ``cracks'' will
inevitably occur if the {\em deficit angle}
\[
\omega_{ABC} \stackrel{\rm def}{=} 2\pi - \vartheta_{ABC}^{\hat F}
- \vartheta_{ABC}^{-\hat E} - \vartheta_{ABC}^{\hat D}
\]
is nonzero.

Now we are ready to write down our ``Euclidean'' formula
\cite[(8)]{TMP4dim1}:
\begin{equation}
V_{\hat D} V_{-\hat E} V_{\hat F} \frac{d\omega_{ABC}}{S_{ABC}} =
V_{\hat A} V_{-\hat B} V_{\hat C} \frac{d\omega_{DEF}}{S_{ABC}}.
\label{eq_z}
\end{equation}
It means the following. First, vertices $A,\ldots,F$ are put
in~$\mathbb R^4$ equipped with Euclidean metric, with the
condition that all six 4-simplices with their vertices in
$A,\ldots,F$ be nondegenerate. This supplies us with edge lengths,
two-dimensional areas, and four-dimensional volumes. Then we vary
the lengths {\em infinitesimally\/} but otherwise {\em
arbitrarily}. This yields an infinitesimal deficit
angle~$d\omega_{ABC}$, as well as a similar infinitesimal
angle~$d\omega_{DEF}$ at the face~$DEF$ around which the
4-simplices $\hat A$, $-\hat B$ and $\hat C$ are situated.

The proof of formula (\ref{eq_z}) can be found
in~\cite[Section~2]{TMP4dim1}. We now emphasize the following: all
quantities entering the l.h.s.\ of~(\ref{eq_z}) belong to the
4-simplices $\hat D$, $\hat E$ and~$\hat F$\footnote{When we are
not preoccupied with being pedantic about the orientation, we drop
the minus sign at $\hat E$ and~$\hat B$.},
 while all quantities entering the r.h.s.\ of~(\ref{eq_z}) belong
 to the 4-simplices $\hat A$, $\hat B$ and~$\hat C$. In this sense, (\ref{eq_z})
 is a relation corresponding to the move~$3\to 3$.

The further way leading to invariants of four-dimensional PL
manifolds goes through the use of acyclic complexes built from
vector spaces consisting of differentials of geometric
quantities~\cite{TMP4dim2,TMP4dim3}. In the present paper, we do
not consider these questions.

\section[The 3 -> 3 relation: affine volume-preserving geometry]{The $\boldsymbol{3\to 3}$
relation: affine volume-preserving geometry}

Our work on invariants based on differential relations between
geometric values attached to the elements of a manifold
triangulation started not from 4-, but from 3-manifolds, and at
first we were using only Euclidean
geometry~\cite{arx_f-la_1,Korepanov_JNMP}. The development of this
subject proceeded, on the one hand, towards four dimensions, as
described in the previous section, and on the other hand ---
towards using other (not Euclidean) geometries. For example,
spherical geometry can be successfully used here, as demonstrated
by Y.~Taylor and C.~Woodward~\cite{TW}. It turned out also that
affine volume-preserving geometry in~$\mathbb R^{n-1}$ can be
used, where $n$ is the dimension of the manifold, which was shown
in papers~\cite{SL2,TMP2004} for $n=3$. In this case, we were able
to introduce some analogues of dihedral angles and deficit angles,
so that the latter measured, so to say, the size of the cracks
appearing when one tries to place a certain geometric
configuration into $\mathbb R^{n-1}=\mathbb R^2$. This geometric
configuration can be imagined as follows: take three tetrahedra
$ABDE$, $BCDE$ and $CADE$ in~$\mathbb R^3$ having a common
edge~$DE$, then project this picture onto~$\mathbb R^2$. Then,
there is a quadratic relation between the areas of projections of
two-dimensional faces. If we, however, allow these areas to be
arbitrary, then our ``deficit angles'' measure to what degree this
quadratic relation is violated. This works also for more than
three tetrahedra situated in a~similar way around a common edge.

In this light, the aim of the present paper looks natural --- to
construct an algebraic relation corresponding to a move $3\to 3$
on a triangulation of a 4-manifold, and based on affine
volume-preserving geometry in~$\mathbb R^3$. As in the previous
section, we will deal with vertices $A,\ldots,F$ and 4-simplices
$\hat A,\ldots,\hat F$. As we take the geometry of~$\mathbb R^3$
for the source of our algebraic relation, we keep in mind two
pictures at once: the four-dimensional combinatorial-topological
picture and the geometrical picture, where {\em three\/}
coordinates are ascribed to each of vertices $A,\ldots,F$ and one
can thus speak of volumes of 3-simplices, i.e.\ tetrahedra. These
volumes will play the key r\^ole below. We will allow ourselves to
denote them by the same letter~$V$ as the 4-volumes in
Section~\ref{sec_recall} (as we are not going to consider any
4-volumes from now on).

Although what follows is a generalization of the result of
paper~\cite{SL2}, using also some ideas from \cite{TMP2004}, our
further exposition will be self-contained. Still, the analogy of
our case $n=4$ with the case $n=3$ of those papers will be one
more guiding thread for us.

There are $6!/(4!\,2!)=15$ tetrahedra with vertices in points
$A,\ldots,F$. Their volumes are not independent. First, there are
{\em linear\/} relations among them: the overall volume of the
boundary of any 4-simplex is zero, e.g.,
\begin{equation}
V_{BCDE}-V_{ACDE}+V_{ABDE}-V_{ABCE}+V_{ABCD}=0. \label{eq_a}
\end{equation}
Certainly, we are always speaking of oriented volumes:
$V_{ABCD}=-V_{BACD}$, etc.

There are five independent relations of type~(\ref{eq_a}).
Besides, there must exist three more relations, which can be seen
from the following counting of parameters (which is no problem to
be made rigorous): 6 points $A,\ldots,F$ have 18 coordinates. From
these, one must subtract the number of parameters in the group of
affine volume-preserving motions of~$\mathbb R^3$, namely 8
parameters corresponding to the subgroup ${\rm SL}(3,\mathbb R)$
and 3 translations. Now everything is all right, the number of
volumes minus the number of relations imposed on them is equal to
the number of coordinates minus the numbers of motions:
\begin{equation}
15-5-3=18-8-3. \label{eq_b}
\end{equation}

We will obtain the three additional (nonlinear) relations among
volumes in such way which will lead us to the natural analogues of
dihedral and deficit angles for our non-metric geometry. Consider
the cluster of three initial 4-simplices $\hat F$, $-\hat E$,
$\hat D$ grouped around their common two-dimensional face~$ABC$.
First, we write down the following easily derived relation which
expresses the vector $\overrightarrow{CE}$ (of course, it lies
in~$\mathbb R^3$) belonging to the 3-face $ABCE$ of
4-simplex~$\hat F$ in terms of the vector $\overrightarrow{CD}$
belonging to the 3-face $ABCD$ of~$\hat F$, the 3-volumes of faces
of~$\hat F$, and two vectors $\overrightarrow{CA}$
and~$\overrightarrow{CB}$ belonging to~$ABC$:
\begin{equation}
\overrightarrow{CE}=\frac{-V_{BCDE}\overrightarrow{CA}+V_{ACDE}\overrightarrow{CB}
+ V_{ABCE}\overrightarrow{CD}}{V_{ABCD}}. \label{eq_c}
\end{equation}
We say briefly that we have expressed $\overrightarrow{CE}$
through $\overrightarrow{CD}$ from 4-simplex~$\hat F$. Similarly,
we express $\overrightarrow{CF}$ through $\overrightarrow{CE}$
from $\hat D$ and $\overrightarrow{CD}$ through
$\overrightarrow{CF}$ from~$\hat E$. Taking the composition of
these expressions, we can say that we go around the face~$ABC$ and
thus express $\overrightarrow{CD}$ through itself. To be exact,
the ``new'' $\overrightarrow{CD}$ comes out as a linear
combination of the ``old'' $\overrightarrow{CD}$ and ``fixed''
vectors $\overrightarrow{CA}$ and~$\overrightarrow{CB}$:
\[
\overrightarrow{CD}_{\rm new} = \overrightarrow{CD}_{\rm old} +
\omega_{ABC}^{(1)}V_{ABCD} \overrightarrow{CA} +
\omega_{ABC}^{(2)}V_{ABCD}\overrightarrow{CB},
\]
where
\begin{equation}
\omega_{ABC}^{(1)} = -\frac{V_{BCDE}}{V_{ABCD}V_{ABCE}} -
\frac{V_{BCEF}}{V_{ABCE}V_{ABCF}}+\frac{V_{BCDF}}{V_{ABCF}V_{ABCD}}
\label{eq_d}
\end{equation}
and
\begin{equation}
\omega_{ABC}^{(2)} = \frac{V_{ACDE}}{V_{ABCD}V_{ABCE}} +
\frac{V_{ACEF}}{V_{ABCE}V_{ABCF}}-\frac{V_{ACDF}}{V_{ABCF}V_{ABCD}}.
\label{eq_e}
\end{equation}

The two nonlinear relations
\begin{equation}
\omega_{ABC}^{(1)} = \omega_{ABC}^{(2)} = 0 \label{eq_f}
\end{equation}
obviously guarantee that the cluster of $\hat F$, $\hat D$
and~$\hat E$, with given 3-volumes of their 3-faces (satisfying,
of course, the conditions of type~(\ref{eq_a})), can be placed
in~$\mathbb R^3$ without contradictions. With nonzero $\omega$'s,
this cannot be done, although each 4-simplex individually still
can be put in~$\mathbb R^3$. We note that each of the quantities
$\omega_{ABC}^{(1)}$ and $\omega_{ABC}^{(2)}$ is a sum of three
expressions of the same kind made of values belonging to
4-simplices $\hat F, \hat D, \hat E$ respectively. We can call
them components of two-component ``dihedral angles'' at the
face~$ABC$ in the mentioned simplices: for instance, the
``dihedral angle'' in~$\hat F$ is
$\Bigl(-V_{BCDE}/(V_{ABCD}V_{ABCE}),\;\allowbreak
V_{ACDE}/(V_{ABCD}V_{ABCE})\Bigr)$. The quantities
$\omega_{ABC}^{(1)}$ and~$\omega_{ABC}^{(2)}$ can be called
components of the ``deficit angle'', or ``discrete curvature''
around~$ABC$.

There must also be a third nonlinear relation between
three-dimensional volumes. As we have shown that it is not
necessary for putting in~$\mathbb R^3$ the 4-simplices $\hat D$,
$\hat E$ and~$\hat F$, it must involve the volumes of 3-faces
absent from these 4-simplices. So, now we consider in the same way
4-simplices $\hat A$, $\hat B$ and~$\hat C$, of which we think as
grouped around the face~$DEF$. We can just change the letters in
formulas (\ref{eq_d}) and (\ref{eq_e}) as follows:
$A\leftrightarrow D$, $B\leftrightarrow E$, $C\leftrightarrow F$,
which means exactly coming from the initial 4-simplices $\hat F$,
$\hat D$ and~$\hat E$ to the simplices $\hat C$, $\hat A$
and~$\hat B$ obtained as a result of the move $3\to 3$. If we now
write out the conditions
\begin{equation}
\omega_{DEF}^{(1)} = \omega_{DEF}^{(2)} = 0 \label{eq_g}
\end{equation}
obtained from (\ref{eq_f}) under this change, we can check (e.g.,
on a computer) that (\ref{eq_f}) and (\ref{eq_g}) together give
exactly three independent conditions.

We now express the three-dimensional volumes through
values~$\lambda$ attached to two-dimen\-sional faces, according to
the following pattern:
\begin{equation}
V_{ABCD}=\lambda_{BCD}-\lambda_{ACD}+\lambda_{ABD}-\lambda_{ABC}.
\label{eq_h}
\end{equation}
In this way we, of course, guarantee that the five linear
relations among volumes hold automa\-tically. In doing so, we act in
analogy with papers~\cite{SL2,TMP2004}, where we had similar
quantities~$\lambda$ attached to {\em edges\/} of triangulation.
Although we do not permute the subscripts of~$\lambda$ in this
paper, it is
 natural to consider $\lambda$ as totally antisymmetric in its indices.

Now we can guess how the algebraic relation corresponding to move
$3\to 3$ can look. In the left-hand side, which corresponds to
simplices $\hat D$, $\hat E$ and~$\hat F$, we expect to find
quantities belonging to the faces specific for that set of
simplices: 2-face $ABC$ and 3-faces $ABCD$, $ABCE$ and~$ABCF$,
while in the right-hand side we expect quantities belonging to the
2-face $DEF$ and 3-faces $ADEF$, $BDEF$ and~$CDEF$. There is no
problem about 3-faces: natural quantities belonging to them are
their volumes. As for the 2-faces, to each of them belong two
components of~$\omega$ and one~$\lambda$, and, as
papers~\cite{SL2,TMP2004} suggest, we must make something like a
derivative $\partial\omega/\partial\lambda$. The right answer
(justified finally by a computer) is to take, for the face $ABC$,
a~Jacobi determinant $J_{ABC,i}$ of partial derivatives of
$\omega_{ABC}^{(1)}$ and~$\omega_{ABC}^{(2)}$ with respect to
$\lambda_{ABC}$ and one more~$\lambda$, call it $\lambda_i$, where
$i$ denotes a two-dimensional face coinciding with neither $ABC$
nor~$DEF$:
\begin{equation}
J_{ABC,i} = \left| \begin{array}{cc} \partial
\omega_{ABC}^{(1)}/\partial
\lambda_{ABC} & \partial \omega_{ABC}^{(1)}/\partial \lambda_{i} \\[.7ex]
\partial \omega_{ABC}^{(2)}/\partial \lambda_{ABC} & \partial \omega_{ABC}^{(2)}/\partial \lambda_{i}
\end{array} \right| .
\label{eq_k}
\end{equation}
All derivatives are taken in the ``flat'' point where all
$\omega$'s are zero and thus all picture is genuinely in~$\mathbb
R^3$.

If we compose also the Jacobi determinant $J_{DEF,i}$ of partial
derivatives of $\omega_{DEF}^{(1)}$ and~$\omega_{DEF}^{(2)}$ with
respect to $\lambda_{DEF}$ and the same $\lambda_i$, then it is
not hard to see that the ratio $J_{ABC,i}/J_{DEF,i}$ does not in
fact depend on the face~$i$. Indeed, if $j$ is another face $\ne
ABC,DEF$, then, according to a version of the theorem on implicit
function derivative,
\begin{equation}
\frac{J_{ABC,i}}{J_{ABC,j}} = -\frac{\partial \lambda_j}{\partial
\lambda_i}, \label{eq_i}
\end{equation}
where the partial derivative in the right-hand side is taken under
the condition that all $\omega$'s are zero. This condition means,
as we explained already, that values~$\lambda$ are such that the
points $A,\ldots,F$ can be placed in~$\mathbb R^3$ in such way
that the tetrahedron volumes coincide with values obtained from
formulas of type~(\ref{eq_h}); and in~(\ref{eq_i}) the three
quantities $\lambda_j$, $\lambda_{ABC}$ and $\lambda_{DEF}$ are
regarded as functions of other $\lambda$'s, including~$\lambda_i$.
It is clear, on the other hand, that $J_{DEF,i}/J_{DEF,j}$ equals
the same right-hand side of~(\ref{eq_i}), so, indeed, the ratio
$J_{ABC,i}/J_{DEF,i}$ is independent of~$i$.

What remains now is to calculate this ratio on a computer and
verify our conjecture that it is expressed in a nice way through
three-dimensional volumes of faces in which our clusters $\hat D$,
$\hat E$, $\hat F$ and $\hat A$, $\hat B$, $\hat C$ differ from
one another. It turns out that these volumes appear in the formula
raised in the second power. Here is this formula --- our formula
corresponding to the move $3\to 3$, proved by a symbolic
calculation using Maple:
\begin{equation}
(V_{ABCD}V_{ABCE}V_{ABCF})^2 J_{ABC,i} =
(V_{ADEF}V_{BDEF}V_{CDEF})^2 J_{DEF,i}. \label{eq_j}
\end{equation}

\section{Concluding remarks}

The result of this note --- formula~(\ref{eq_j}) --- relies upon
computer-assisted calculations which were being made in the belief
that a four-dimensional analogue of formula~(15) from~\cite{SL2}
must exist. Thus, a careful investigation of algebraic structures
arising here is still to be done.

We emphasize once more that we consider the constructing of
various new relations, such as~(\ref{eq_j}), and investigating the
related algebraic structures as a way to discovering new quantum
topological field theories, in particular, for manifolds of higher
$(>3)$ dimensions.

\subsection*{Acknowledgements}

I would like to use this occasion to thank the organizers of the
conference ``Symmetry in Nonlinear Mathematical Physics'' in Kiev
in June, 2005, for their excellent conference and inviting me to
write this paper. My special thanks to S.~Podobedov, Member of the
Parliament of Ukraine, for his warm hospitality during my stay in
Kiev.

This paper was written with a partial financial support from
Russian Foundation for Basic Research, Grant no.~04-01-96010.

\LastPageEnding


\newpage


\begin{thebibliography}{99}
\footnotesize

\bibitem{Lickorish}
Lickorish W.B.R., Simplicial moves on complexes and manifolds,
{\it Geom. Topol. Monographs}, 1999, Vol.~2, Proceedings of the
Kirbyfest, 299--320; math.GT/9911256.

\bibitem{leftsymmetric}
Diatta A., Medina A., Classical Yang--Baxter equation and left
invariant affine geometry on Lie groups. {\it Manuscripta
mathematica}, 2004, V.114, N~4, 477--486; math.DG/0203198.

\bibitem{arx_f-la_1}
Korepanov I.G., A formula with volumes of five tetrahedra and
discrete curvature, nlin.SI/0003001.

\bibitem{arx_f-la_2}
Korepanov I.G., A formula with hypervolumes of six 4-simplices and
two discrete curvatures, nlin.SI/0003024.

\bibitem{Korepanov_JNMP}
Korepanov I.G., Invariants of PL manifolds from metrized
simplicial complexes. Three-dimen\-sional case, {\it J.~Nonlinear
Math. Phys.}, 2001, V.8, N~2, 196--210; math.GT/0009225.

\bibitem{SL2}
Korepanov I.G., Martyushev E.V., A classical solution of the
pentagon equation related to the group $SL(2)$, {\it Theor. Math.
Phys.}, 2001, V.129, N~1, 1320--1324.

\bibitem{Ponzano-Regge}
Ponzano G., Regge T., Semi-classical limit of Racah coefficients,
in ``Spectroscopic and Group Theoretical Methods in Physics'',
Editor F. Bloch, North-Holland, 1968, 1--58.

\bibitem{Roberts}
Roberts J., Classical $6j$-symbols and the tetrahedron, {\it Geom.
Topol.}, 1999, V.3, 21--66; math-ph/9812013.

\bibitem{TW}
Taylor Y., Woodward C., Spherical tetrahedra and invariants of
3-manifolds, math.GT/0406228.

\bibitem{TMP2004}
Korepanov I.G., SL(2)-Solution of the pentagon equation and
invariants of three-dimen\-sional manifolds, {\it Theor. Math.
Phys.}, 2004, V.138, N~1, 18--27; math.AT/0304149.

\bibitem{TMP4dim1}
Korepanov I.G., Euclidean 4-simplices and invariants of
four-dimen\-sional manifolds: I.~Moves $3\to 3$, {\it Theor. Math.
Phys.}, 2002, V.131, N~3, 765--774; math.GT/0211165.

\bibitem{TMP4dim2}
Korepanov I.G., Euclidean 4-simplices and invariants of
four-dimen\-sional manifolds: II.~An algebraic complex and moves
$2\leftrightarrow 4$, {\it Theor. Math. Phys.}, 2002, V.133, N~1,
1338--1347; math.GT/0211166.

\bibitem{TMP4dim3}
Korepanov I.G., Euclidean 4-simplices and invariants of
four-dimen\-sional manifolds: III.~Moves $1\leftrightarrow 5$ and
related structures, {\it Theor. Math. Phys.}, 2003, V.135, N~2,
601--613; math.GT/0211167.

\bibitem{Atiyah}
Atiyah M.F., Topological quantum field theory, {\it Publications
Math\'ematiques de l'IH\'ES}, 1988, V.68, 175--186.


\end{thebibliography}
\end{document}